\documentclass[12pt]{article}%
\usepackage{amsmath}
\usepackage{amsfonts}
\usepackage{amssymb}
\usepackage{graphicx}%
\providecommand{\U}[1]{\protect\rule{.1in}{.1in}}
\newtheorem{theorem}{Theorem}[section]

\newtheorem{definition}[theorem]{Definition}

\newtheorem{remark}[theorem]{Remark}

\numberwithin{equation}{section}

\begin{document}

\title{Synthesis of a quantum tree Weyl matrix}
\author{Sergei A. Avdonin{\small $^{\text{1}}$}, Kira V. Khmelnytskaya$^{\text{2}}$,
Vladislav V. Kravchenko{\small $^{\text{3}}$}\\{\small $^{\text{1}}$ Department of Mathematics and Statistics, University of
Alaska, Fairbanks, AK 99775, USA}\\{\small $^{\text{2}}$ Faculty of Engineering, Autonomous University of
Queretaro, }\\{\small Cerro de las Campanas s/n, col. Las Campanas Quer\'{e}taro, Qro. C.P.
76010 M\'{e}xico}\\{\small $^{\text{3}}$ Department of Mathematics, Cinvestav, Campus
Quer\'{e}taro, }\\{\small Libramiento Norponiente \#2000, Fracc. Real de Juriquilla,
Quer\'{e}taro, Qro., 76230 M\'{e}xico}\\{\small e-mail: s.avdonin@alaska.edu, khmel@uaq.edu.mx,
vkravchenko@math.cinvestav.edu.mx}}
\maketitle

\begin{abstract}
A method for successive synthesis of a Weyl matrix (or Dirichlet-to-Neumann
map) of an arbitrary quantum tree is proposed. It allows one, starting from
one boundary edge, to compute the Weyl matrix of a whole quantum graph by
adding on new edges and solving elementary systems of linear algebraic
equations in each step.

\end{abstract}

\section{Introduction}

Quantum graphs or differential equations networks have wide applications in
science and engineering and give rise to challenging problems involving many
areas of modern mathematics, from combinatorics to partial differential
equations and spectral theory. A number of surveys and collections of papers
on quantum graphs appeared last years, including the first books on this topic
by Berkolaiko and Kuchment \cite{BerkolaikoKuchment} and Mugnolo
\cite{Mugnolo}. In the present work we consider tree graphs, that is, finite
connected compact graphs without cycles. The Weyl or Titchmarsh-Weyl matrix of
a quantum graph is one of the key mathematical objects, it naturally appears
in direct and inverse spectral theory, in the control theory of quantum graphs
and in numerous applications. Its importance lies in the fact that the Weyl
matrix (or more precisely, its transposed matrix) is the Dirichlet-to-Neumann
map of the quantum graph. For a fixed value of the spectral parameter, a
vector of arbitrary Dirichlet-type boundary values of a solution multiplied by
the transposed Weyl matrix gives one a vector of the corresponding
Neumann-type boundary values, if only the value of the spectral parameter is
not a Dirichlet eigenvalue of the quantum graph. Moreover, the singularities
of the Weyl matrix, considered as a function of the spectral parameter,
determine the Dirichlet spectrum of the quantum graph. Many papers on inverse
problems for quantum graphs exploit the Weyl matrix or equivalent spectral
data, see, e.g., \cite{Be}, \cite{BrWe}, \cite{Yu}, \cite{BeVa}, \cite{FrYu},
\cite{AvdoninKurasov2008}.

Direct construction of the Weyl matrix for a sufficiently large quantum tree
is quite a challenging problem. It requires solving large systems of equations
involving solutions (and their derivatives) of differential equations on all
edges of the tree. The main result of the present work is a simple procedure
for a progressive synthesis of the Weyl matrix of an arbitrary quantum tree.
Starting from just one leaf edge and adding successively new edges, it allows
one to compute the Weyl matrices for ever larger quantum trees from the Weyl
matrices for the smaller ones. We call this procedure the synthesis of the
Weyl matrix. In a sense, it is inverse with respect to the leaf peeling
method, which was developed in \cite{AvdoninKurasov2008}, see also
\cite{AvdoninZhao2020}, \cite{AKK2023tree}. The proposed synthesis of the Weyl
matrix is based on revealed relations between Weyl solutions of a larger
quantum tree graph $\Omega$ with those of a smaller one $\widetilde{\Omega}$,
obtained from $\Omega$ by cutting out all the leaf edges of an internal vertex.

In Section \ref{Sect Preliminaries} we recall necessary definitions and
express the Weyl solutions in terms of fundamental systems of solutions of the
Sturm--Liouville equation on each edge. Section \ref{Sect Synthesis} presents
the main result of the work, the procedure of the synthesis of a Weyl matrix
of an arbitrary quantum tree graph. Finally, Section \ref{Sect Concl} contains
some concluding remarks.

\section{Preliminaries\label{Sect Preliminaries}}

Let $\Omega$ be a finite connected compact graph without cycles (a tree graph)
consisting of $P$ edges, $e_{1}, \ldots, e_{P},$ and $P+1$ vertices,
$V=\left\{  v_{1},...,v_{P+1}\right\}  .$ The notation $e_{j}\sim v$ means
that the edge $e_{j}$ is incident to the vertex $v$. Every edge $e_{j}$ is
identified with an interval $(0,L_{j})$ of the real line. The boundary
$\Gamma=\left\{  \gamma_{1},\ldots,\gamma_{m}\right\}  $ of $\Omega$ is the
set of all leaves of the graph (the external vertices). The edge adjacent to
some $\gamma_{j}\ $is called a leaf or boundary edge.

A continuous function $u$ defined on the graph $\Omega$ is a $P$-tuple of
functions $u_{j}\in C\left[  0,L_{j}\right]  $ satisfying the continuity
condition at the internal vertices $v$: $u_{i}(v)=u_{j}(v)$ for all
$e_{i},e_{j}\sim v$. Then $u\in C(\Omega)$.

Let $q\in%
\mathcal{L}%
_{1}(\Omega)$ be real valued, and $\lambda$ a complex number. Consider the
Sturm-Liouville equation on $\Omega$:
\begin{equation}
-u^{\prime\prime}(x)+q(x)u(x)=\lambda u(x). \label{Schr}%
\end{equation}
A function $u$ defined on $\Omega$ is said to be a solution of (\ref{Schr}) if
besides (\ref{Schr}) we have that
\begin{equation}
u\in C(\Omega), \label{continuity}%
\end{equation}
and for every internal vertex $v$ the Kirchhoff-Neumann condition is
fulfilled
\begin{equation}
\sum_{e_{j}\sim v}\partial u_{j}(v)=0,\text{\quad for all }v\in V\setminus
\Gamma. \label{KN}%
\end{equation}
Here $u_{j}$ is a restriction of $u$ onto $e_{j}$, $\partial u_{j}(v)$ stands
for the derivative of $u$ at the vertex $v$ taken along the edge $e_{j}$ in
the direction outward the vertex, and the sum is taken over all the edges
incident to the internal vertex $v$.

For simplicity we assume that the considered graph does not have vertices of
degree two, because every such vertex can be regarded as an internal point of
an edge which is a sum of two edges incident at such a vertex, and the
continuity condition (\ref{continuity}) together with the Kirchhoff-Neumann
condition (\ref{KN}) guarantee that any solution of (\ref{Schr}) on the
incident edges keeps satisfying (\ref{Schr}) also on the union of them.

\begin{definition}
A solution $w_{i}$ of (\ref{Schr}) on $\Omega$ is called the \textbf{Weyl
solution} associated with the leaf $\gamma_{i}$ if it satisfies the boundary
conditions
\begin{equation}
w_{i}(\gamma_{i})=1\quad\text{and}\quad w_{i}(\gamma_{j})=0\text{ for all
}j\neq i. \label{WeylBoundary}%
\end{equation}

\end{definition}

If $\lambda$ in (\ref{Schr}) is not a Dirichlet eigenvalue of the quantum tree
$\Omega$, the Weyl solution $w_{i}$ exists and is unique for any
$i=1,\ldots,m$. In particular, since the potential $q$ is real valued, the
Dirichlet spectrum is real and thus the boundary value problem (\ref{Schr}),
(\ref{WeylBoundary}) is uniquely solvable for all $\lambda\notin\mathbb{R}$.

\begin{definition}
The $m\times m$ matrix-function $\mathbf{M}(\lambda)$, $\lambda\notin
\mathbb{R}$, consisting of the elements $\mathbf{M}_{ij}(\lambda)=\partial
w_{i}(\gamma_{j})$, $i,j=1,\ldots,m$ is called the \textbf{Weyl matrix}.
\end{definition}

For a fixed value of $\lambda$, the transposed Weyl matrix represents a
Dirichlet-to-Neumann map of the quantum graph defined by $\Omega$ and $q\in%
\mathcal{L}%
_{1}(\Omega)$. Indeed, if $u$ is a solution of (\ref{Schr}) satisfying the
Dirichlet condition at the boundary vertices $u(\lambda,\gamma)=f(\lambda)$,
then $\partial u(\lambda,\gamma)=\mathbf{M}^{T}(\lambda)f(\lambda)$,
$\lambda\notin\mathbb{R}$.

It is clear that the direct computation of the Weyl matrix, which involves
finding Weyl solutions and computing their derivatives at leaves may be a
difficult task, especially, when $\Omega$ consists of a large number of edges.
The corresponding systems of equations which need to be solved in this case
may be too large, because they should combine information on the solutions and
their derivatives on each edge and at all internal vertices.

The main result of the present work is a simple procedure which allows one to
synthesize the Weyl matrix progressively, by adding edges to smaller graphs
and computing the Weyl matrices for the obtained larger graphs from the Weyl
matrices for the smaller ones. We call this procedure synthesis of the Weyl
matrix. It allows one to compute the Weyl matrix of any quantum tree starting
from one leaf edge and adding successively new edges.

By $\varphi_{i}(\rho,x)$ and $S_{i}(\rho,x)$ we denote the so-called
fundamental solutions of the Sturm--Liouville equation on the edge $e_{j}:$
\begin{equation}
-y^{\prime\prime}(x)+q_{i}(x)y(x)=\rho^{2}y(x),\quad x\in(0,L_{i}),
\label{Schri}%
\end{equation}
satisfying the initial conditions
\[
\varphi_{i}(\rho,0)=1,\quad\varphi_{i}^{\prime}(\rho,0)=0,
\]%
\[
S_{i}(\rho,0)=0,\quad S_{i}^{\prime}(\rho,0)=1.
\]
Here $q_{i}(x)$ is the component of the potential $q(x)$ on the edge $e_{i}$,
and $\rho=\sqrt{\lambda}$, $\operatorname{Im}\rho\geq0$. For a leaf edge
$e_{i}$ it is convenient to identify its leaf $\gamma_{i}$ with the left
endpoint $x=0.$ Then the Weyl solution $w_{i}(\rho,x)$ has the form%
\[
w_{ii}(\rho,x)=\varphi_{i}(\rho,x)+\mathbf{M}_{i,i}(\rho^{2})S_{i}%
(\rho,x)\quad\text{on the adjacent leaf edge }e_{i}%
\]
and
\[
w_{ij}(\rho,x)=\mathbf{M}_{i,j}(\rho^{2})S_{j}(\rho,x)\quad\text{on every
other leaf edge }e_{j},\quad j\neq i.
\]
Hereafter, the notation $w_{ij}(\rho,x)$ means that we consider $j$-th
component of a solution $w_{i}(\rho,x)$, that is, the solution $w_{i}(\rho,x)$
on the edge $e_{j}$.

On internal edges $e_{j}$ we have
\[
w_{ij}(\rho,x)=a_{ij}(\rho)\varphi_{j}(\rho,x)+b_{ij}(\rho)S_{j}(\rho,x),
\]
where the choice of which vertex is identified with zero is arbitrary, and in
general the factors $a_{ij}(\rho)$, $b_{ij}(\rho)$ are unknown.

Since in direct and inverse spectral problems involving the Weyl matrix one
deals with solutions on large ranges of the parameter $\rho$, it is convenient
to use the Neumann series of Bessel functions representations for $\varphi
_{j}(\rho,x)$ and $S_{j}(\rho,x)$, introduced in \cite{KNT} and applied in a
number of direct and inverse problems (see, e.g., \cite{AKK2023},
\cite{AKK2023tree}, \cite{AvdoninKravchenko2022}, \cite{KrBook2020},
\cite{KKC2022Mathematics}). One of the features of these representations is
the existence of estimates for the remainders of the series independent of
$\operatorname{Re}\rho$.

\section{Synthesis of Weyl matrix\label{Sect Synthesis}}

Consider a quantum tree graph $\widetilde{\Omega}$ whose leaves are
$\gamma_{0},\gamma_{1},\ldots,\gamma_{m}$. Assume its Weyl matrix
$\widetilde{\mathbf{M}}(\rho^{2})$ to be known for some value of $\rho$.
Attach a number of edges to the leaf $\gamma_{0}$ (see Fig. 1), so that
$\gamma_{0}$ becomes an internal vertex of a new larger graph $\Omega$, and
$\gamma_{1},\ldots,\gamma_{m},\gamma_{m+1},\ldots,\gamma_{m+m_{1}}$ the leaves
of $\Omega$. Here $m_{1}$ is the number of the new attached edges. Thus,
$\Omega$ is a tree graph obtained from $\widetilde{\Omega}$ by attaching
$m_{1}$ edges to $\gamma_{0}$. For simplicity we call these new edges
$e_{m+1},\ldots,e_{m+m_{1}}$ and denote their respective lengths as
$L_{m+1},\ldots,L_{m+m_{1}}$. We assume that a corresponding potential
$q_{j}\in\mathcal{L}_{1}(0,L_{j})$, $j=m+1,\ldots,m+m_{1},$ is given on each
new edge, and the leaf $\gamma_{j}$ is identified with $x=0.$
%

\begin{figure}
[ptb]
\begin{center}
\includegraphics[
height=2.9915in,
width=4.4211in
]%
{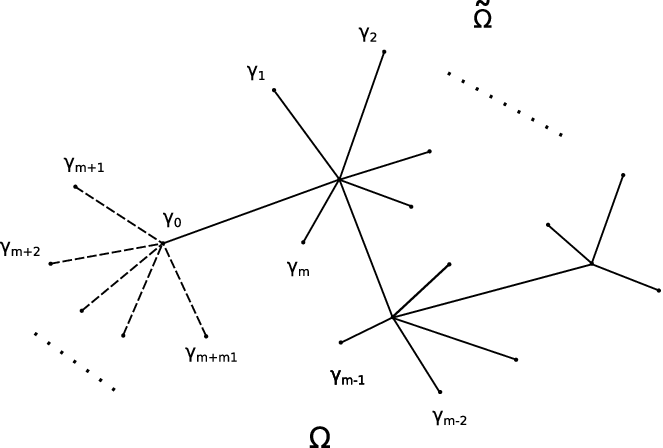}%
\caption{Tree graph $\Omega$ is obtained from a subgraph $\widetilde{\Omega}$
(its edges are presented by solid lines) by attaching to the vertex
$\gamma_{0}$ a number of new edges (dashed lines).}%
\label{Fig1}%
\end{center}
\end{figure}

Our task is to find the Weyl matrix $\mathbf{M}(\rho^{2})$ of the quantum
graph $\Omega$. The idea is to construct the Weyl solutions $w_{i}(\rho,x)$ of
$\Omega$ in terms of the Weyl solutions $\widetilde{w}_{j}(\rho,x)$ of
$\widetilde{\Omega}$. We start with the Weyl solutions $w_{m+j}(\rho,x)$
associated with the new leaves $\gamma_{m+j}$, $j=1,\ldots,m_{1}$.

Let us look for $w_{m+j}(\rho,x)$ in the form%
\begin{equation}
w_{m+j}(\rho,x)=c_{j}(\rho)\widetilde{w}_{0}(\rho,x)\quad\text{on }%
\widetilde{\Omega}. \label{w_m+j}%
\end{equation}
That is, on the subgraph $\widetilde{\Omega}$ the Weyl solution $w_{m+j}%
(\rho,x)$ coincides with $\widetilde{w}_{0}(\rho,x)$ up to a multiplicative
constant $c_{j}(\rho)$. In this case, it automatically satisfies the
homogeneous Dirichlet condition at $\gamma_{1},\ldots,\gamma_{m},$ and we
still need to satisfy the conditions $w_{m+j}(\rho,\gamma_{m+j})=1$ and
$w_{m+j}(\rho,\gamma_{m+i})=0$ for $i\neq j$.

From (\ref{w_m+j}) we have
\begin{equation}
\varphi_{m+j}(\rho,L_{m+j})+\mathbf{M}_{m+j,m+j}(\rho^{2})S_{m+j}(\rho
,L_{m+j})=c_{j}(\rho) \label{continuity1}%
\end{equation}
and
\begin{equation}
\varphi_{m+j}^{\prime}(\rho,L_{m+j})+\sum_{k=1}^{m_{1}}\mathbf{M}%
_{m+j,m+k}(\rho^{2})S_{m+k}^{\prime}(\rho,L_{m+k})=c_{j}(\rho)\widetilde
{\mathbf{M}}_{0,0}(\rho^{2}), \label{KN m+j}%
\end{equation}
where $\widetilde{\mathbf{M}}_{0,0}(\rho^{2})$ is an element of the Weyl
matrix $\widetilde{\mathbf{M}}(\rho^{2})$: $\widetilde{\mathbf{M}}_{0,0}%
(\rho^{2})=\partial\widetilde{w}_{0}(\gamma_{0})$.

Additionally, the continuity condition gives us the equalities%
\begin{equation}
\mathbf{M}_{m+j,m+k}(\rho^{2})S_{m+k}(\rho,L_{m+k})=c_{j}(\rho),\quad
k=1,\ldots,m_{1}\text{ and }k\neq j. \label{continuity2}%
\end{equation}

Thus, for each $j$, from (\ref{continuity1}), (\ref{KN m+j}) and
(\ref{continuity2}) we have $m_{1}+1$ equations for the $m_{1}+1$ unknowns%
\[
\left\{  \mathbf{M}_{m+j,m+k}(\rho^{2}),\,k=1,\ldots,m_{1};\text{ }c_{j}%
(\rho)\right\}  .
\]
Note that in the linear algebraic system (\ref{continuity1}), (\ref{KN m+j}),
(\ref{continuity2}) the magnitudes $\varphi_{m+j}(\rho,L_{m+j})$,
$S_{m+k}(\rho,L_{m+k})$, $\varphi_{m+j}^{\prime}(\rho,L_{m+j})$,
$S_{m+k}^{\prime}(\rho,L_{m+k})$ are known, since all the potentials
$q_{m+k}(x)$, $k=1,\ldots,m_{1}$ are known.

Thus, from (\ref{continuity1}), (\ref{KN m+j}) and (\ref{continuity2}) we find
$\mathbf{M}_{m+j,m+k}(\rho^{2}),\,k=1,\ldots,m_{1}$ and $c_{j}(\rho)$.

From (\ref{w_m+j}) we obtain additionally,
\begin{equation}
\mathbf{M}_{m+j,i}(\rho^{2})=c_{j}(\rho)\widetilde{\mathbf{M}}_{0,i}(\rho
^{2})\text{\quad for }i=1,\ldots,m, \label{M_m+j,i}%
\end{equation}
and thus we have already completed the rows $m+1,\ldots,m+m_{1}$ of the Weyl
matrix $\mathbf{M}(\rho^{2})$.

Now, choose an $i\in\left\{  1,\ldots,m\right\}  $. We look for $w_{i}%
(\rho,x)$ such that on $\widetilde{\Omega}$ the equality be valid%
\[
w_{i}(\rho,x)=\widetilde{w}_{i}(\rho,x)+\alpha_{i}(\rho)\widetilde{w}_{0}%
(\rho,x),
\]
where $\alpha_{i}(\rho)$ is a constant. This is a natural choice, because
\[
w_{i}(\rho,\gamma_{i})=1\quad\text{and}\quad w_{i}(\rho,\gamma_{j}%
)=0\text{\quad for }j=1,\ldots,m\text{ and }j\neq i.
\]
Moreover, we have
\[
w_{i}(\rho,\gamma_{0})=\alpha_{i}(\rho)
\]
and
\[
\partial w_{i}(\rho,\gamma_{0})=\widetilde{\mathbf{M}}_{i,0}(\rho^{2}%
)+\alpha_{i}(\rho)\widetilde{\mathbf{M}}_{0,0}(\rho^{2}).
\]
Thus, for all $j=m+1,\ldots,m+m_{1}$ we have%
\begin{equation}
w_{ij}(\rho,\gamma_{0})=\alpha_{i}(\rho) \label{w_ij}%
\end{equation}
and%
\begin{equation}
\sum_{j=m+1}^{m+m_{1}}\partial w_{ij}(\rho,\gamma_{0})=-\widetilde{\mathbf{M}%
}_{i,0}(\rho^{2})-\alpha_{i}(\rho)\widetilde{\mathbf{M}}_{0,0}(\rho^{2}).
\label{KN_ij}%
\end{equation}
Equality (\ref{w_ij}) can be written as%
\begin{equation}
\mathbf{M}_{i,j}(\rho^{2})S_{j}(\rho,L_{j})=\alpha_{i}(\rho),\quad
j=m+1,\ldots,m+m_{1}, \label{system2_1}%
\end{equation}
while (\ref{KN_ij}) takes the form%
\begin{equation}
\sum_{j=m+1}^{m+m_{1}}\mathbf{M}_{i,j}(\rho^{2})S_{j}^{\prime}(\rho
,L_{j})=\widetilde{\mathbf{M}}_{i,0}(\rho^{2})+\alpha_{i}(\rho)\widetilde
{\mathbf{M}}_{0,0}(\rho^{2}). \label{system2_2}%
\end{equation}
For every $i\in\left\{  1,\ldots,m\right\}  $, equations (\ref{system2_1}) and
(\ref{system2_2}) give us $m_{1}+1$ equations for the $m_{1}+1$ unknowns%
\[
\left\{  \mathbf{M}_{i,j}(\rho^{2}),\,j=m+1,\ldots,m+m_{1};\text{ }\alpha
_{i}(\rho)\right\}  .
\]
The magnitudes $S_{j}(\rho,L_{j})$, $S_{j}^{\prime}(\rho,L_{j})$,
$\widetilde{\mathbf{M}}_{i,0}(\rho^{2})$, $\widetilde{\mathbf{M}}_{0,0}%
(\rho^{2})$ in (\ref{system2_1}) and (\ref{system2_2}) are already known.

Finally, to obtain $\mathbf{M}_{i,j}(\rho^{2})$ for $i,j=1,\ldots,m$ we
observe that
\[
\partial w_{i}(\rho,\gamma_{j})=\partial\widetilde{w}_{i}(\rho,\gamma
_{j})+\alpha_{i}(\rho)\partial\widetilde{w}_{0}(\rho,\gamma_{j})
\]
and thus%
\begin{equation}
\mathbf{M}_{i,j}(\rho^{2})=\widetilde{\mathbf{M}}_{i,j}(\rho^{2})+\alpha
_{i}(\rho)\widetilde{\mathbf{M}}_{0,j}(\rho^{2})\quad\text{for }%
i,j=1,\ldots,m\text{.} \label{M_ij}%
\end{equation}

Let us summarize the procedure of construction of the Weyl matrix
$\mathbf{M}(\rho^{2})$ of the quantum tree $\Omega$ from the Weyl matrix
$\widetilde{\mathbf{M}}(\rho^{2})$ of the quantum subtree $\widetilde{\Omega}$.

1) For each $j=1,\ldots,m_{1}$ solve the $\left(  m_{1}+1\right)
\times\left(  m_{1}+1\right)  $-system of linear algebraic equations
(\ref{continuity1})-(\ref{continuity2}) to find the constant $c_{j}(\rho)$ and
the Weyl matrix entries $\mathbf{M}_{m+j,m+k}(\rho^{2}),\,k=1,\ldots,m_{1}$.
Compute the entries $\mathbf{M}_{m+j,i}(\rho^{2}),$ $i=1,\ldots,m$ from
(\ref{M_m+j,i}). Thus, the rows from $m+1$ to $m+m_{1}$ of the Weyl matrix
$\mathbf{M}(\rho^{2})$ are computed.

2) For each $i=1,\ldots,m$ solve the $\left(  m_{1}+1\right)  \times\left(
m_{1}+1\right)  $-system of linear algebraic equations (\ref{system2_1}),
(\ref{system2_2}) to find the constant $\alpha_{i}(\rho)$ and the Weyl matrix
entries $\mathbf{M}_{i,j}(\rho^{2})$, $j=m+1,\ldots,m+m_{1}$. Compute the
entries $\mathbf{M}_{i,j}(\rho^{2}),$ $j=1,\ldots,m$ from (\ref{M_ij}). This
completes the construction of the Weyl matrix $\mathbf{M}(\rho^{2})$.

Finally, let us consider the situation when $\widetilde{\Omega}$ is just a
single segment with the vertices $\gamma_{0}$ and $\gamma_{1}$. Since the
potential on this segment is supposed to be given, we may assume that both
corresponding Weyl solutions on such $\widetilde{\Omega}$ are known:
$\widetilde{w}_{0}(\rho,x)$ and $\widetilde{w}_{1}(\rho,x)$, that satisfy the
boundary conditions
\[
\widetilde{w}_{0}(\rho,\gamma_{0})=1,\quad\widetilde{w}_{0}(\rho,\gamma_{1})=0
\]
and
\[
\widetilde{w}_{1}(\rho,\gamma_{0})=0,\quad\widetilde{w}_{1}(\rho,\gamma
_{1})=1.
\]
Then the entries of the $2\times2$ - Weyl matrix $\widetilde{\mathbf{M}}%
(\rho^{2})$ have the form%
\[
\widetilde{\mathbf{M}}_{0,0}(\rho^{2})=\partial\widetilde{w}_{0}(\rho
,\gamma_{0}),\quad\widetilde{\mathbf{M}}_{0,1}(\rho^{2})=\partial\widetilde
{w}_{0}(\rho,\gamma_{1}),
\]%
\[
\widetilde{\mathbf{M}}_{1,0}(\rho^{2})=\partial\widetilde{w}_{1}(\rho
,\gamma_{0}),\quad\widetilde{\mathbf{M}}_{1,1}(\rho^{2})=\partial\widetilde
{w}_{1}(\rho,\gamma_{1}).
\]
Now, as a first step for synthesis of the Weyl matrix of a quantum tree we
attach a number of edges to the vertex $\gamma_{0}$, that gives us a quantum
star graph. The procedure described above gives us its Weyl matrix.
Subsequently, attaching new edges to the leaves and applying the above
procedure leads to the computation of the Weyl matrix of an ever larger
quantum tree. Thus, starting with one edge we synthesize the Weyl matrix of
the whole quantum tree.

\begin{remark}
The requirement on the potential $q$ to be real valued is not essential. The
proposed procedure for the synthesis of the Weyl matrix is applicable to
complex valued potentials without modifications. The only constraint that must
be checked at each step is that $\rho\in\mathbb{C}$ does not belong to the
Dirichlet spectrum of either graph $\widetilde{\Omega}$ or $\Omega$.  
\end{remark}

\section{Conclusions\label{Sect Concl}}

A procedure of the synthesis of a Weyl matrix of an arbitrary quantum tree
graph is developed. It allows one to compute the Weyl matrix for a large
quantum tree successively, by adding new edges and computing the Weyl matrices
for ever larger subgraphs. In each step relatively small systems of linear
algebraic equations are solved. Due to the fact that the transposed Weyl
matrix is the Dirichlet-to-Neumann map of the quantum tree graph, the
synthesis procedure will find applications in solving a variety of boundary
value and control problems on quantum tree graphs.

\textbf{Funding }The research of Sergei Avdonin was supported in part by the
National Science Foundation, grant DMS 1909869, and by Moscow Center for
Fundamental and Applied Mathematics. The research of Vladislav Kravchenko was
supported by CONACYT, Mexico, via the project 284470.

\textbf{Data availability} The data that support the findings of this study
are available upon reasonable request.

\textbf{Declarations}

\textbf{Conflict of interest} The authors declare no competing interests.

\end{document}